

\magnification\magstephalf
\vsize=24.0truecm
\baselineskip13pt






\def\hatt{\widehat}

\def\tilda{\widetilde}
\def\eps{\varepsilon}

\def\half{\hbox{$1\over2$}}

\def\quart{\hbox{$1\over4$}}
\def\sixth{\hbox{$1\over6$}}
\def\arr{\rightarrow}

\def\RR{\mathord{I\kern-.3em R}}
\def\PP{\mathord{I\kern-.3em P}}
\def\NN{\mathord{I\kern-.3em N}}
\def\ZZ{\mathord{I\kern-.3em Z}} 

\def\E{{\rm E}}

\def\subsection{\medskip}

\font\bigbf=cmbx12

\font\csc=cmcsc10
\font\cyr=wncyr10
 at 10truept 
\font\smallrm=cmr8

\def\today{\number\day \space \ifcase\month\or
January\or February\or March\or April\or May\or June\or 
July\or August\or September\or October\or November\or December\fi  
\space \number\year}


   
\def\ref#1{{\noindent\hangafter=1\hangindent=20pt
  #1\smallskip}}          

\def\quotationone{\smallrm Where there is a Will}
\def\quotationtwo{\smallrm There is a Won't}
\def\hskipdistanceleft{\hskip-3.5pt}
\def\hskipdistanceright{\hskip-2.0pt}
\footline={{
\ifodd\count0
        {\hskipdistanceleft\quotationone\phantom{\smallrm\today}
                \hfil{\rm\the\pageno}\hfil
         \phantom{\quotationone}{\smallrm\today}\hskipdistanceright}
        \else 
        {\hskipdistanceleft\quotationtwo\phantom{\today}
                \hfil{\rm\the\pageno}\hfil
         \phantom{\quotationtwo}{\smallrm\today}\hskipdistanceright}
        \fi}}

         
\def\cstok#1{\leavevmode\thinspace\hbox{\vrule\vtop{\vbox{\hrule\kern1pt
        \hbox{\vphantom{\tt/}\thinspace{\tt#1}\thinspace}}
        \kern1pt\hrule}\vrule}\thinspace} 
\def\square{\cstok{\phantom{$\cdot$}}} 


\def\fermat#1{\setbox0=\vtop{\hsize4.00pc
        \smallrm\raggedright\noindent\baselineskip9pt
        \rightskip=0.5pc plus 1.5pc #1}\leavevmode
        \vadjust{\dimen0=\dp0
        \kern-\ht0\hbox{\kern-4.00pc\box0}\kern-\dimen0}}

\def\hsizeplusepsilon{14.25truecm} 
\def\fermatright#1{\setbox0=\vtop{\hsize4.00pc
        \smallrm\raggedright\noindent\baselineskip9pt
        \rightskip=0.5pc plus 1.5pc #1}\leavevmode
        \vadjust{\dimen0=\dp0
        \kern-\ht0\hbox{\kern\hsizeplusepsilon\box0}\kern-\dimen0}}



\def\quotationone{\smallrm Hjort and Fenstad}
\def\quotationtwo{\smallrm Second order asymptotics}
\def\today{September 1994} 

\def\j3ord{\rm\u{\cyr i}}

\centerline{\bigbf Second order asymptotics for the number of times}
\centerline{\bigbf an estimator is more than $\eps$ from its target value} 

\medskip 
\centerline{\bf Nils Lid Hjort and Grete Fenstad, University of Oslo} 


\smallskip 
{{\smallskip\narrower\baselineskip12pt\noindent
{\csc Abstract.} Suppose $\{\hatt\theta_n\colon n\ge1\}$ 
is a strongly consistent sequence of estimators 
for a parameter $\theta$, where 
$\hatt\theta_n$ is based on the first $n$ observations.
Consider $Q_\eps$, the number of times $|\hatt\theta_n-\theta|\ge\eps$.
In another paper (Hjort and Fenstad, 1992) 
we have shown that $\eps^2Q_\eps$ has a limit distribution
as $\eps\rightarrow0$, depending only on $\sigma$,
the standard deviation of the limit distribution for 
$\sqrt{n}(\hatt\theta_n-\theta)$,
under natural regularity conditions. 
The present paper investigates some second order asymptotics for 
differences between $Q_\eps$-variables. 
The limit of $\E(Q_{1,\eps}-Q_{2,\eps})$ is calculated in cases where 
$\E Q_{1,\eps}/\E Q_{2,\eps}$ goes to 1, 
leading to a notion of `asymptotic relative deficiency'
in cases where the asymptotic relative efficiency is 1.
This is used to distinguish between 
competing estimators with identical limit distributions. 
Thus using denominator $n-{1\over3}$ 
in the familiar formula for estimating a
normal variance is better than both $n$ and $n-1$ and indeed 
all other choices, for example,
in the sense of leading to the smallest possible 
expected number of $\eps$-errors. 
Results of this type are found in a selection of 
familiar estimation problems,
using limit results for expected differences,
and are compared to corresponding asymptotic relative deficiency 
analysis in the sense of Hodges and Lehmann. 
Some second order distributional results are reached as well.
It is shown how $\eps$ times a $Q_\eps$-difference 
tends to a variable which is related to some 
exponential distributions associated with Brownian motion,
and that have recently been investigated by 
Hjort and Khasminskii (1993). 

\smallskip\noindent
{\csc Key words:} {\sl 
asymptotic relative deficiency,
asymptotic relative efficiency,
comparison of estimators, 
loss function,
second order optimality,
total relative time for Brownian motion, 
the number of $\eps$-misses }
\smallskip}}

\bigskip

{\bf 1. Introduction.}
Suppose $X_1,X_2,\ldots$ is a sequence of independent 
observations from some distribution $F$, and let $\theta=\theta(F)$
be some parameter of interest. 
Assume that $\hatt\theta_n$ 
is based on the first $n$ data points, 
and that this estimator sequence is 
strongly consistent, i.e.~$\hatt\theta_n\rightarrow\theta$ almost surely. 
Then $Q_\eps$, the number of cases where $|\hatt\theta_n-\theta|\ge\eps$, 
is finite almost surely. Under natural regularity conditions, 
which include existence of a normal $(0,\sigma^2)$ limit 
for $\sqrt{n}(\hatt\theta_n-\theta)$, 
Hjort and Fenstad (1992, Section 7) have shown that 
$$\eps^2Q_\eps\rightarrow_d Q
	={\rm Leb}\{s\ge0\colon |W(s)|\ge s/\sigma\} 
	=\int_0^\infty 
	I\{|W(s)|\ge s/\sigma\}\,{\rm d}s, \eqno(1.1)$$
in which ${\rm Leb}$ is Lebesgue measure on the {half}line,
$W(.)$ is a Brownian motion process, and $I\{...\}$ denotes
an indicator function. 
They also show that $\eps^2\E Q_\eps\rightarrow \E Q$,
which simply is equal to $\sigma^2$. 
Suppose $\sqrt{n}(\hatt\theta_{n,j}-\theta)\rightarrow_d N\{0,\sigma_j^2\}$
for two competing estimator sequences,
and let $Q_{j,\eps}$ be their accompanying number of $\eps$-misses
variables. 
Then a natural asymptotic measure of relative efficiency is 
$${\rm a.r.e.}=\lim_{\eps\rightarrow0}
	{\E Q_{1,\eps}\over \E Q_{2,\eps}}
	={\sigma_1^2\over \sigma_2^2}. \eqno(1.2)$$
This is also the traditional formula for a.r.e.~based on 
calculations of asymptotic ratios $n_1/n_2$ in which
$n_1$ and $n_2$ are the sample sizes needed by methods 1 and 2 to
achieve the same level of mean squared error precision. 
Hjort and Fenstad (1992) also give generalisations of (1.1) and (1.2)
to multi-dimensional parameters and general distance measures,
and reach limit results for $Q_\eps$-variables 
in more involved problems like nonparametric estimation
of distribution functions and densities. 

The $Q_\eps$-based criterion (1.2) is only a `first order measure'
and cannot distinguish between
estimator sequences with identical limit distributions,
in which case typically $\eps^2(Q_{1,\eps}-Q_{2,\eps})\rightarrow0$
and $Q_{1,\eps}/Q_{2,\eps}\rightarrow1$ in probability. 
Our aim in this paper is to develop some `second order' theory for 
{\it differences} of $Q_\eps$'s, which then can be applied to 
single out the `best' estimator in a class of estimators with
the same limit distribution. 
Being `best' here means having 
the smallest possible expected number of $\eps$-errors, 
in the limit as $\eps$ goes to zero. 
In particular we intend to solve comparison problems for 
unbiased minimum variance (UMV) estimators versus 
maximum likelihood (ML) estimators versus Bayes estimators 
in some familiar cases where they disagree. 

This is in a spirit similar to 
Hodges and Lehmann (1970) who worked with 
a measure of asymptotic relative deficiency, 
expressed as the limit of sample size difference $n_1-n_2$ 
in cases where $n_1/n_2\rightarrow1$; 
see also Lehmann (1983, Chapter 5.2). 
In Section 2 formulae are found for a $Q_\eps$-based 
$${\rm a.r.d.}=\lim_{\eps\rightarrow0}
		\E(Q_{1,\eps}-Q_{2,\eps}), \eqno(1.3)$$ 
in the simple case of estimating a mean parameter, 
using estimators of the Bayesian variety 
${n\over n+c}\bar X_n+{c\over n+c}d$. 
In this a.r.d.~formula the underlying skewness 
of the $X_i$'s enters in a natural way,
whereas it does not in Hodges and Lehmann's calculations. 
The a.r.d.~formula makes it possible to show that 
$n-{1\over3}$ is the superior choice of 
denominator in the familiar $\sum_{i=1}^n(X_i-\bar X_n)^2/(n+c)$ 
formula for a normal variance, for example. 
This and other mean-parameter related 
examples are treated in Section 3. 
A couple of technically more demanding 
problems are included in Section 4. 
We are able to show that $(\bar X_n)^2+\hatt\sigma_n^2/n$,
as an estimator for the squared normal mean,   
can be expected to exhibit fewer $\eps$-errors than all 
other estimators of the form 
$(\bar X_n)^2+d\hatt\sigma_n^2/n$, for example.

Comparing estimation methods in terms of $\E Q_\eps$ 
can be phrased in decision theoretic terms, 
where the underlying loss function 
is equal to the number of $\eps$-errors 
for the full sequence of estimates. 
This point is briefly discussed in Section 5, 
along with some Bayesian considerations. 

The theory and applications presented in Sections 2--5 
are only concerned with {\it expected values} of $Q_\eps$-variables.
This is reflected in the tools used, 
namely Edgeworth expansions and further Taylor type 
approximations to probabilities. 
Section 6 presents a brief addendum and 
is concerned with second order {\it distributional} 
aspects of $Q_{1,\eps}-Q_{2,\eps}$ differences.
In the structurally simplest case of estimating a mean parameter
it turns out that 
$$\eps^2(Q_{1,\eps}-Q_{2,\eps})\rightarrow_p0
	{\rm \ and\ } Q_{1,\eps}/Q_{2,\eps}\rightarrow_p1 
	\quad {\rm while} \quad 
	\eps(Q_{1,\eps}-Q_{2,\eps})\rightarrow_d A-B, $$
say, where $A$ and $B$ are certain total relative time 
variables associated with Brownian motion, with 
distributions that are exponential or mixtures of exponentials
and unit point masses at zero; 
see Hjort and Khasminskii (1993) for a separate account on these. 

Our paper is mainly reporting on an investigation into 
one particular way of distinguishing between first order
equivalent sequences. 
We do not claim that the a.r.d.~criterion (1.3) 
on which our second order comparisons are based is 
statistically more natural than other criteria,
but it is interesting that it can be computed at all,
and to see how it fares compared to other second order criteria, 
of which several others exist.
In addition to the a.r.d.-calculations of Hodges and Lehmann
and ours of the present paper, 
Rao (1962) and others have worked with asymptotic sufficiency,
Pfanzagl (1973) and others with coverage probabilities and 
median unbiasedness, and Ghosh and Subramanyam (1974) 
and others with expansions for mean squared errors. 
A good source for further information is the
discussion of the paper by Berkson (1980). 
Let us finally point out that techniques of the present paper
can be used to establish expansion results for 
coverage probabilities of the type
$p_n={\rm Pr}\{|\hatt\theta_n-\theta|\le k/\sqrt{n}\}$.
These results, which in a way are simpler than those of Sections 2--5,
could then be used to compare different estimators 
with the same limit distribution. 
Such comparisons would however depend on the value of~$k$. 

\bigskip
{\bf 2. General results for estimating a mean.}
Suppose $X_1,X_2,X_3,\ldots$ are i.i.d.~with 
$\E X_i=\xi$, ${\rm Var}\,X_i=\sigma^2$, and 
skewness $\E(X_i-\xi)^3/\sigma^3=\gamma$. 
The most useful and directly informative results are 
those of Propositions 1 and 2 below, but for technical reasons 
it is convenient to start with the following lemma. 

\smallskip
{\csc Lemma.} 
{\sl Let $a$ be positive, 
and consider $Q_\eps(c)$, the number of times, among $n\ge a/\eps^2$, 
where $|{n\over n+c}\bar X_n-\xi|\ge\eps$.
Suppose that $X_i$ has finite fourth order moment
and that its distribution is non-lattice, 
and write $\phi(.)$ for the standard normal density. Then}
$$\eqalign{
\lambda_a(c)
	&=\lim_{\eps\rightarrow0}\E\{Q_\eps(c)-Q_\eps(0)\} \cr
	&=2\Bigl[{\xi^2\over\sigma^2}c^2
	-\Bigl(2-{2\gamma\over3}{\xi\over\sigma}\Bigr)c\Bigr]
	\int_{\sqrt{a}/\sigma}^\infty \phi(y)\,{\rm d}y
	-{2\gamma\over3}{\xi\over\sigma}{\sqrt{a}\over\sigma}
		\phi(\sqrt{a}/\sigma)\,c. \cr} \eqno(2.1)$$

{\csc Proof:} 
We choose to work with $Q_\eps(c)$ in the form 
$$Q_\eps(c)=\sum_{n\ge am}I\Bigl\{\Big|{n\over n+c}\bar X_n-\xi\Big|\ge
	{1\over \sqrt{m}}\Bigr\}, \eqno(2.2)$$
writing $m=1/\eps^2$. Note that $Q_\eps(c)$ is finite a.s.~by the 
strong law of large numbers. 
Consider $T_n=\sqrt{n}(\bar X_n-\xi)/\sigma$. 
Inserting $\bar X_n=\xi+\sigma T_n/\sqrt{n}$ 
the indicator function term in (2.2) can be written 
$I\{T_n\le l(c){\rm \ or\ }T_n\ge r(c)\}$, where 
$$l(c)=-{1\over\sigma}\sqrt{n\over m}
	+{1\over \sqrt{n}}{c\xi\over \sigma}
	-{1\over \sqrt{nm}}{c\over \sigma}, 
	\quad 
  r(c)={1\over\sigma}\sqrt{n\over m}
	+{1\over \sqrt{n}}{c\xi\over \sigma}
	+{1\over \sqrt{nm}}{c\over \sigma}. \eqno(2.3)$$
Letting $s=n/m$ we have 
$l(0)=-u$, 
$l(c)=-u+b_1/\sqrt{n}-b_2/n$, 
$r(0)=u$, 
$r(c)=u+b_1/\sqrt{n}+b_2/n$, 
in which 
$u=\sqrt{s}/\sigma$, $b_1=c\xi/\sigma$, and $b_2=c\sqrt{s}/\sigma$. 
Since 
$$Q_\eps(c)-Q_\eps(0)=\sum_{n\ge am}
	\Bigl[I\bigl\{T_n\le l(c){\rm \ or\ }T_n\ge r(c)\bigr\}
	     -I\bigl\{T_n\le l(0){\rm \ or\ }T_n\ge r(0)\bigr\}\Bigr] 
				\eqno(2.4)$$
its expected value can be written
$$E\{Q_\eps(c)-Q_\eps(0)\}
	=\sum_{n\ge am}\Bigl[\bigl\{G_n(l(c))-G_n(l(0))\bigr\}
		-\bigl\{G_n(r(c))-G_n(r(0))\bigr\}\Bigr], \eqno(2.5)$$
in terms of the cumulative distribution function $G_n$ of $T_n$. 

Under the non-lattice assumption 
there is a Cram\'er--Edgeworth expansion for $G_n$ of the form 
$$G_n(t)=\Phi(t)-{1\over6}{\gamma\over \sqrt{n}}A(t)\phi(t)
		+{1\over n}R(t)\phi(t)+O(n^{-3/2}), \eqno(2.6)$$
where $A(t)=t^2-1$ and $R(t)$ is a certain polynomial of degree five, 
see for example Barndorff-Nielsen and Cox (1989, Ch.~4).
In view of (2.5) we need to approximate terms of the type
$G_n(t+\delta)-G_n(t)$ where $\delta$ is of the order of $1/\sqrt{n}$.
Using Taylor expansions one finds that 
$$G_n(t+\delta)-G_n(t)=\phi(t)\delta-{1\over2}t\phi(t)\delta^2
	+{1\over6}{\gamma\over\sqrt{n}}B(t)\phi(t)\delta+O(n^{-3/2}),$$ 
in which $B(t)\phi(t)$ is the derivative of $-A(t)\Phi(t)$, 
i.e.~$B(t)=t^3-3t$, and the $O$-term is uniform in $t$. This leads to 
$$\eqalign{
G_n(-u+b_1/\sqrt{n}-b_2/n)-G_n(-u)
	&=\phi(u)b_1/\sqrt{n}+\phi(u)\{-b_2+\half ub_1^2
		-\sixth\gamma b_1B(u)\}/n, \cr
G_n(u+b_1/\sqrt{n}+b_2/n)-G_n(u)
	&=\phi(u)b_1/\sqrt{n}+\phi(u)\{b_2-\half ua^2
		+\sixth\gamma b_1B(u)\}/n. \cr}$$
The summand in (2.5) can accordingly be expressed,
apart from $O(n^{-3/2})$ terms, as
$${1\over n}\phi(u)\Bigl\{-2b_2+ub_1^2-{1\over3}\gamma b_1B(u)\Bigr\}
	={1\over m}\phi\Bigl({\sqrt{s}\over\sigma}\Bigr)
	 {1\over s}\Bigl\{-2{\sqrt{s}\over \sigma}c
	 +{\xi^2\sqrt{s}\over \sigma^3}c^2-{\gamma\over3}{\xi\over \sigma} 
	 B\Bigl({\sqrt{s}\over \sigma}\Bigr)c\Bigr\}. $$
Since ${1\over m}\sum_{n/m\ge a}h(n/m)$ is a Riemannian approximation
to the integral $\int_a^\infty h(s)\,{\rm d}s$, and converges 
to this limit as $m\rightarrow\infty$, we find in the end
that $\E\{Q_\eps(c)-Q_\eps(0)\}$ converges to 
$$\eqalign{
\lambda_a(c)
&=\int_a^\infty\Bigl[{\xi^2\over\sigma^3}{1\over\sqrt{s}}c^2
	-{2\over\sigma}{1\over\sqrt{s}}c
	-{\gamma\over3}{\xi\over\sigma}{1\over s}
	  B\Bigl({\sqrt{s}\over \sigma}\Bigr)c\Bigr]
		\,\phi\Bigl({\sqrt{s}\over \sigma}\Bigr)\,{\rm d}s \cr
&=2\int_{\sqrt{a}/\sigma}^\infty
	\Bigl[{\xi^2\over\sigma^2}c^2-2c
		-{\gamma\over3}{\xi\over\sigma}{B(u)\over u}c\Bigr]
		\,\phi(u)\,{\rm d}u. \cr} $$ 
Some further analysis, using $B(u)/u=u^2-3$, 
finally proves (2.1). \square 

\smallskip
In the applications presented below we will study this 
and similar limits as functions of $c$, 
to find the estimator sequence that 
can be expected to make fewest $\eps$-errors. 
The optimal value of $c$ as computed from (2.1) 
will depend upon the somewhat arbitrary value of $a$, however. 
Write in general $Q_{\eps,m}(c)$ for the number of $\eps$-errors 
committed by ${n\over n+c}\bar X_n$ among $n\ge m$ cases. 
The limit result (2.1) for $\lambda_a(c)$ 
relates to $Q_{\eps,a/\eps^2}(c)$,
and it is natural to let $a$ tend to zero.
The limit function $\lambda_0(c)$ obtained 
by letting $a\rightarrow0$ in (2.1)  
is not quite the limit of the expected difference 
between $Q_{\eps,1}(c)$ and $Q_{\eps,1}(0)$, since the remainder term 
$O(\sum_{n\ge am}n^{-3/2})$ does not go to zero for $a=1/m$.
The following holds, however. 

\smallskip 
{\csc Proposition 1.} 
{\sl Let $Q_{\eps}(c)$ be the number 
of $\eps$-misses for ${n\over n+c}\bar X_n$, 
counted this time among $n\ge a(\eps)/\eps^2$, where 
$a(\eps)\rightarrow0$ while $a(\eps)/\eps^2\rightarrow\infty$
(as with $a(\eps)=\eps$, for example). Then} 
$$\lambda_0(c)
	=\lim_{\eps\rightarrow0}\E\bigl\{Q_{\eps,a(\eps)/\eps^2}(c)
	-Q_{\eps,a(\eps)/\eps^2}(0)\bigr\}
	={\xi^2\over \sigma^2}\,c^2
     	-2\Bigl(1-{\gamma\over3}{\xi\over \sigma}\Bigr)\,c. \eqno(2.7)$$ 

\smallskip
Sometimes it is reasonable to balance the objective $\bar X_n$ against
some prior value. Bayesians often use estimators 
of the type $\hatt\xi_n(c,d)={n\over n+c}\bar X_n+{c\over n+c}\xi_0$,
see Section 5. Consider therefore the slight extension 
$$Q_\eps(c,d)=\sum_{n\ge a/\eps^2}I\Bigl\{\Big|
	{n\over n+c}\bar X_n+{c\over n+c}d-\xi\Big|\ge\eps\Bigr\}. $$
The appropriate generalisation of (2.3) is found to be 
$$l(c,d)=-{1\over\sigma}\sqrt{n\over m}
	+{1\over \sqrt{n}}{c(\xi-d)\over \sigma}
	-{1\over \sqrt{nm}}{c\over \sigma}, 
	\quad 
  r(c,d)={1\over\sigma}\sqrt{n\over m}
	+{1\over \sqrt{n}}{c(\xi-d)\over \sigma}
	+{1\over \sqrt{nm}}{c\over \sigma}, $$
and by repeating previous arguments mutatis mutandis one proves 
the following generalisation of (2.7). 

\smallskip
{\csc Proposition 2.} 
{\sl Let as above $\eps$ and $a(\eps)=\eps$ tend to zero. Then} 
$$\lambda_0(c,d)
	=\lim_{\eps\rightarrow0}\E\bigl\{Q_\eps(c,d)-Q_\eps(0,0)\bigr\}
	={(\xi-d)^2\over \sigma^2}\,c^2
     	-2\Bigl(1-{\gamma\over3}{\xi-d \over \sigma}\Bigr)\,c. \eqno(2.8)$$ 

\smallskip 
A further generalisation of (2.7) and (2.8) is given at the end of 4C below.
We think of these as a.r.d.-formulae, see (1.3) and the remarks made there.
We also remark that there is an alternative way of proving 
(2.1), (2.7), (2.8), 
via distributional analysis following the Brownian motion approximations
briefly discussed in Section 6. Details pertaining to this
alternative proof are in the technical report 
Hjort and Fenstad (1993) available from the authors.
\eject 

\smallskip
{\csc Remark.} 
One may also compute the Hodges and Lehmann asymptotic relative deficiency
in this situation, defined as the limit of $n_0(c,d)-n_0$,
where $n_0(c,d)$ is the sample size needed to achieve 
$\E\{\hatt\xi_n(c,d)-\xi\}^2=\E\{\hatt\xi_{n_0}(0,0)-\xi\}^2$, say; 
see Lehmann (1983, Section 5.2) for discussion. 
Then $n_0(c,d)/n_0(0,0)\rightarrow1$ for all $(c,d)$, 
but one can prove 
$${\rm a.r.d.}_{\rm hl}(c,d)
	=\lim_{n_0\rightarrow\infty}\{n_0(c,d)-n_0(0,0)\}
	={(\xi-d)^2\over \sigma^2}c^2-2c. \eqno(2.9)$$
The a.r.d.~formula (2.8) differs from this in that the skewness 
$\gamma$ also enters, in a natural way. \square 

\bigskip
{\bf 3. Some applications.}

\smallskip
{\sl 3A. Normal mean.}
Let $X_i$ be normal $(\theta,1)$. 
Consider the estimator sequence 
$\hatt\theta_n(c,d)={n\over n+c}\bar X_n+{c\over n+c}d$,
and let $Q_\eps(c,d)$ be the number $\eps$-misses. 
This fits into the general scheme 
of Section 2 with $\xi=\theta$, $\sigma=1$, and $\gamma=0$. Hence 
$$\lambda_0(c,d)
	=\lim_{\eps\rightarrow0}\E_\theta
		\bigl\{Q_\eps(c,d)-Q_\eps(0,0)\bigr\}
	=(\theta-d)^2c^2-2c.$$
If this is averaged over $\theta$ w.r.t.~some distribution 
with `prior mean' $\E\theta=\theta_0$ and `prior variance'
${\rm Var}\,\theta=\tau^2$, then one finds 
$\{\tau^2+(\theta_0-d)^2\}c^2-2c$, 
which is minimised when $d$ is chosen as $\theta_0$ 
and $c$ is chosen as $1/\tau^2$. 
This provides fresh and independent motivation for using 
$$\theta_n^*={n\over n+1/\tau^2}\bar X_n
	+{1/\tau^2\over n+1/\tau^2}\theta_0 \eqno(3.1)$$
in the presence of such prior knowledge, 
and agrees with both familiar Bayesian calculations in the normal model
and the so-called credibility formula of actuarial statistics. 
Observe that $\theta_n^*$ can expect to make fewer $\eps$-errors
than $\bar X_n$ does if the true $\theta$ parameter is within
$\sqrt{2}\tau$ of $\theta_0$. 
Note also that only $\theta_0$ and $\tau^2$ matter 
regarding the choice of prior weight function.  
Further discussion about risks and average risks is in Section 5. 


\smallskip
{\sl 3B. Exponential mean.} 
Let $X_i\sim {\rm Exp}(1/\theta)$, and consider $Q_\eps(c)$, 
the number of times
$|\hatt\theta_n(c)/\theta-1|\ge\eps$,
where $\hatt\theta_n(c)={n\over n+c}\bar X_n$. 
This can be written $|{n\over n+c}\bar Y_n-1|\ge\eps$, 
where the $Y_i$'s come from $e^{-y}$. This is as in Section 2
with $\xi=1$, $\sigma=1$, $\gamma=2$. The limit of interest is 
$c^2-{2\over3}c$, with best value $c_0={1\over 3}$. 
The ML solution with $c=0$ can be expected to make $1\over9$ 
more $\eps$-errors [sic] while the best estimator under squared error loss,
which uses $c=1$, can be expected to make $4\over9$ more $\eps$-errors. 

\smallskip
{\sl 3C. Normal variance.} 
Next consider $\hatt\sigma_N^2(c)=\sum_{i=1}^N(Y_i-\bar Y_N)^2/(N-1+c)$
for estimating the variance $\sigma^2$ 
based on data $Y_1,\ldots,Y_N$ that are normal $(\mu,\sigma^2)$.  
Study the number of times $|\hatt\sigma_N^2(c)/\sigma^2-1|\ge\eps$.
This can be written 
$|{n\over n+c}(\chi^2_{n}/n)-1|\ge\eps$, writing $n$ for $N-1$,
the degrees of freedom. This is once more as in Section 2, 
this time with $\chi^2_1$ variables playing the r\^ole of $X_i$'s. 
These have $\xi=1$, $\sigma=\sqrt{2}$ and $\gamma=2\sqrt{2}$. 
The limit to study is ${1\over2}c^2-{2\over3}c$, 
with minimum occurring for $c_0={2\over 3}$. Hence 
$$\hatt\sigma_N^2={1\over N-(1/3)}\sum_{i=1}^N(Y_i-\bar Y_N)^2 \eqno(3.2)$$
can expect to make the fewest $\eps$-errors! 
See also Section 4C below. 


\smallskip
{\sl 3D. Binomial probability.}
Let $Y_n$ be binomial $(n,p)$, and let $Q_\eps(c,d)$ count $\eps$-misses
for $(Y_n+cd)/(n+c)$. Note that $Y_n$ is the sum of $X_i$'s with 
mean $p$, variance $pq$, and skewness $\gamma=(q-p)/(pq)^{1/2}$, 
where $q=1-p$. Hence the risk difference limit in (2.8) becomes 
$${\rm risk}(c,d)={(p-d)^2\over pq}c^2-2c
	-{2\over3}{(p-q)(p-d)\over pq}c. $$
After studying the case with `prior guess' $d=\half$ one 
can show that  
$p_n^*=(Y_n+2/3)/(n+4/3)$ 
is the second order minimax sequence of estimators, 
in the sense of minimising the maximum possible risk.  
It can, seemingly, expect to make 2.667 fewer $\eps$-errors than $Y_n/n$ does,
regardless of the underlying $p\,$! 

There is an obstacle here, however, in that the distribution 
of $X_i$ is lattice, the Cram\'er-condition does not hold,
and fine-tuned corrections are necessary for a formula like 
(2.6) to hold. Techniques from Kolassa and McCullagh (1990) 
are likely to provide in the end a formula like (2.8) appropriate
for the lattice case, but we have not pursued this. 
A possible trick to avoid these difficulties associated 
with exact a.r.d.~calculations here is to spread 
the probability mass $p$ for $X_i$ 
uniformly on some tiny $[1-\eta,1+\eta]$.
Now results of Section~2 are in force, 
giving the best sequence of estimates for $p$.
Letting $\eta\rightarrow0$ in the end singles out the 
$(Y_n+2/3)/(n+4/3)$ sequence as the best. 
Note finally that the ${\rm Beta}(2/3,2/3)$ 
is a least favourable prior for $p$, 
in the sense of the calculations above.

\bigskip
{\bf 4. Further examples.} 
The aim of this section is to establish results in a couple of 
technically more demanding estimation problems. 
In each case somewhat strenuous 
modifications of the arguments used in Section 2 have to be devised
in order to establish the limit of expected $Q_\eps$-differences. 

\smallskip
{\sl 4A. The squared mean in the normal model: known variance.} 
Let $X_i$ be normal $(\xi,\sigma^2)$, with $\sigma$ known, 
and suppose $\theta=\xi^2$ is to be estimated. 
The ML solution is $(\bar X_n)^2$, which could
overestimate; the UMV solution is $(\bar X_n)^2-\sigma^2/n$.
Study therefore $\hatt\theta_n(d)=(\bar X_n)^2-d\sigma^2/n$
and its corresponding $Q_\eps(d)$, the number of $\eps$-misses
among $n\ge a(\eps)/\eps^2$, with $a(\eps)$ going to zero as in (2.7). 
We intend to show that 
$$\lambda_0(d)
	=\lim_{\eps\rightarrow0}\E\bigl\{Q_\eps(d)-Q_\eps(0)\bigr\}
	= (\quart d^2+\half d){\sigma^2\over \xi^2}. \eqno(4.1)$$

Some alterations to the programme of Section 2 are necessary. 
The indicator function summand in $Q_\eps(d)$ is 1 if
$|(\xi+T_n\sigma_n/\sqrt{n})^2-d\sigma^2/n-\xi^2|\ge 1/\sqrt{m}$,
in which $T_n=\sqrt{n}(\bar X_n-\xi)/\sigma$. This can be written
$$T_n^2+2(\xi/\sigma)\sqrt{n}T_n\ge d+\sqrt{n}\sqrt{s}/\sigma^2
	\quad {\rm or} \quad
  T_n^2+2(\xi/\sigma^2)\sqrt{n}T_n\le d-\sqrt{n}\sqrt{s}/\sigma^2,$$
where $s=n/m$ again. Let $b(d)<r(d)$ be the two roots of the 
first inequality and $c(d)<l(d)$ the two roots of the second inequality.
Assume now that $\xi$ is positive. Then 
$b(d)<c(d)<l(d)<0<r(d)$, and the indicator function term is 
$$I\{T_n\le b(d){\rm\ or\ }T_n\ge r(d){\rm\ or\ }c(d)\le T_n\le l(d)\}.$$
It turns out that both $b(d)$ and $c(d)$ are of size 
about $-2(\xi/\sigma)\sqrt{n}$,
which is too far out on the left to be of significance, and we can
concentrate on $I\{T_n\le l(d){\rm \ or\ }T_n\ge r(d)\}$. 
We need to compute the limit of 
$$\E\{Q_\eps(d)-Q_\eps(0)\}
	\doteq\sum_{n\ge am}\Bigl[\bigl\{\Phi(l(d))-\Phi(l(0))\bigr\}
		-\bigl\{\Phi(r(d))-\Phi(r(0))\bigr\}\Bigr], \eqno(4.2)$$
cf.~(2.4) and (2.5).
A determined mind finds 
$$\eqalign{
l(d)&=-\sqrt{n}{\xi\over\sigma}+\sqrt{n}{\xi\over\sigma}
	\Bigl(1-{1\over \sqrt{n}}{\sqrt{s}\over \xi^2}
	+{1\over n}{d\sigma^2\over \xi^2}\Bigr)^{1/2} \cr
    &\doteq-{1\over2}{\sqrt{s}\over \xi\sigma}
	+{1\over\sqrt{n}}
         \Bigl({1\over2}{d\sigma\over\xi}-{1\over8}{s\over\xi^3\sigma}\Bigr)
	+{1\over n}
	 \Bigl({1\over4}d\sqrt{s}{\sigma\over \xi^3}
		-{1\over 16}{s^{3/2}\over \xi^5\sigma}\Bigr) \cr
    &=-u+\delta_1(d)/\sqrt{n}+\delta_2(d)/n, \cr}$$
$$\eqalign{       
r(d)&=-\sqrt{n}{\xi\over \sigma}+\sqrt{n}{\xi\over\sigma}
	\Bigl(1+{1\over \sqrt{n}}{\sqrt{s}\over \xi^2}
	+{1\over n}{d\sigma^2\over \xi^2}\Bigr)^{1/2} \cr
    &\doteq{1\over2}{\sqrt{s}\over \xi\sigma}
	+{1\over\sqrt{n}}
         \Bigl({1\over2}{d\sigma\over\xi}-{1\over8}{s\over\xi^3\sigma}\Bigr)
	+{1\over n}
	 \Bigl(-{1\over4}d\sqrt{s}{\sigma\over \xi^3}
		+{1\over 16}{s^{3/2}\over \xi^5\sigma}\Bigr) \cr
    &=u+\delta_1(d)/\sqrt{n}-\delta_2(d)/n. \cr}$$ 
One can now compute Taylor approximations to the terms of (4.2).
The result, apart from $O(n^{-3/2})$-terms, is
$$\eqalign{
2\phi(u)\{&\delta_2(d)-\delta_2(0)\}/n
	+\{\half\phi'(-u)-\half\phi'(u)\}\{\delta_1(d)^2-\delta_1(0)^2\}/n \cr
	&={1\over ms}\phi\Bigl({1\over2}{\sqrt{s}\over \xi\sigma}\Bigr)
	 \Bigl[{1\over2}d\sqrt{s}{\sigma\over \xi^3}
	       +\Bigl({1\over4}{d^2\sigma^2\over \xi^2}
		-{1\over8}{ds\over \xi^4}\Bigr)
		 {1\over2}{\sqrt{s}\over \xi\sigma}\Bigr]. \cr}$$
Further care and attention to details lead to (4.1).

When $\xi<0$ then similar arguments give the same answer 
$(\quart d^2+\half d)\,\sigma^2/\xi^2$. 
We may conclude that the choice $d_0=-1$ is best: 
Within the class of estimators under consideration, 
the estimator $(\bar X_n)^2+\sigma^2/n$ has the 
second order optimality property of making the fewest $\eps$-errors, 
in expectation, in the limit as $\eps\rightarrow0$.
The ML solution $(\bar X_n)^2$ can expect to make $\quart\sigma^2/\xi^2$ 
more $\eps$-errors while the UMV $(\bar X_n)^2-\sigma^2/n$ 
can expect to make $\sigma^2/\xi^2$ more $\eps$-errors. 

A natural extension is to include all estimators of the type 
$\hatt\theta_n(c,d)=({n\over n+c}\bar X_n)^2-d\sigma^2/n$. 
Analysis is provided in Hjort and Fenstad (1993),
giving the best possible values of $c$ and $d$ 
in terms of `prior averages' of $\xi^2/\sigma^2$ and $\sigma^2/\xi^2$;
see also Section 5 below. Included in Hjort and Fenstad (1993) 
is also analysis of the special case $\xi=0$, where the asymptotics 
work differently. 
%


\smallskip
{\sl 4B. The squared mean in the normal model: unknown variance.}
The extension to unknown $\sigma$ is important because this
parameter usually {\it is} unknown, of course, 
but also because the techniques given in a moment serve to illustrate 
how similar problems can be solved 
in other two- and multi-parameter situations.  

The natural class of estimators to study is 
$\theta_n^*(d)=(\bar X_n)^2-d\hatt\sigma_n^2/n$,
in which we take $\hatt\sigma_n^2=\sum_{i=1}^n(X_i-\bar X_n)^2/(n-1)$ 
to be the unbiased version. Again $d=0$ gives the ML solution while
$d=1$ gives the UMV estimator. 
Write $\hatt\sigma_n^2=\sigma^2Z_n$, 
where $Z_n$ is independent of $\bar X_n$ and distributed as 
$\chi^2_{n-1}/(n-1)$, so that in fact $\theta_n^*(d)=\hatt\theta_n(dZ_n)$,
in the notation of 4A. The arguments there show that $Q^*_\eps(d)$,
the number of cases where $|\theta_n^*(d)-\xi^2|\ge\eps$, 
essentially counts cases where $T_n\le l(dZ_n)$ or $T_n\ge r(dZ_n)$. 
The analogue of (4.2) becomes 
$$
  \sum_{n\ge am}\Bigl\{\bigl[{\rm Pr}\{T_n\le l(dZ_n)\}
	-{\rm Pr}\{T_n\le l(0)\}\bigr]
	-\bigl[{\rm Pr}\{T_n\le r(dZ_n)\}
	-{\rm Pr}\{T_n\le r(0)\}\bigr]\Bigr\}.$$
The simplest way to compute the limit of this sum, as $m=1/\eps^2$ grows,
is to condition on the value of $Z_n$, use results of 4A, and then
integrate over the distribution $g_n(z_n)$ for $Z_n$.
Carrying through this gives 
$$\sum_{n\ge am}\int_0^\infty
	{1\over ms}\phi\Bigl({1\over2}{\sqrt{s}\over \xi\sigma}\Bigr)
	\Bigl[{1\over2}d\sqrt{s}z_n{\sigma\over \xi^3}
	+{1\over2}{\sqrt{s}\over \xi\sigma}\Bigl(
		{1\over4}d^2z_n^2{\sigma^2\over \xi^2}
		-{1\over8}dz_n{s\over \xi^4}\Bigr)\Bigr]
			\,g_n(z_n)\,{\rm d}z_n.$$
But since $EZ_n=1$ and $EZ_n^2=1+2/(n-1)$ the limit 
becomes in the end equal to the previous answer 
$(\quart d^2+\half d)\,\sigma^2/\xi^2$. We may conclude that 
$(\bar X_n)^2+\hatt\sigma^2_n/n$ is best. 

It is interesting to compare this result to 
the corresponding one using the Hodges and Lehmann deficiency.
Calculations with the mean squared error of $\theta_n^*(d)$ 
show that ${\rm a.r.d.}_{\rm hl}$, the limit of the sample
size difference $n_0(d)-n_0(0)$, becomes  
$(\quart d^2-\half d)\,\sigma^2/\xi^2$, and the best 
estimator with this criterion is the UMV solution 
$(\bar X_n)^2-\hatt\sigma^2_n/n$. 

\smallskip
{\sl 4C. Standard deviation in the normal model.} 
Let us next present a variant of Example 3C. 
Instead of counting instances of 
$|\hatt\sigma_N^2(c)/\sigma^2-1|\ge\eps$, 
penalise errors using the natural scale for $\sigma$ and study 
$$Q_\eps(c)=\sum_{N\ge a/\eps^2}
	I\{|\hatt\sigma_N(c)/\sigma-1|\ge\eps\}. \eqno(4.3)$$
The finer aspects of the second order asymptotics machinery 
turn out to give a different answer than in 3C for the best value of $c$. 

To find the best value this time we provide a modest but useful 
generalisation of (2.7) and (2.8). Suppose $h(\xi)$ 
is a smooth increasing transformation, and consider $Q^*_\eps(c,d)$, 
the number of times $|h(\hatt\xi_n(c,d))-h(\xi)|\ge\eps$,
among $n\ge 1/\eps$ (see (2.7)),  
with $\hatt\xi_n(c,d)=(n\bar X_n+cd)/(n+c)$ as in Proposition 2
of Section 2. One can go through previous arguments, 
supplement them with further Taylor analysis, and prove 
$$\lambda_0^*(c,d)
	=\lim_{\eps\rightarrow0}\E\{Q^*_\eps(c,d)-Q^*_\eps(0,0)\}
	={(\xi-d)^2\over \sigma^2}c^2
	 +\Bigl\{-2+{2\gamma\over3}{\xi-d\over\sigma}
		+{k_2\over k_1^2}(\xi-d)\Bigr\}\,c, \eqno(4.4)$$
in which $k_1=(h^{-1})'(h(\xi))$ and $k_2=(h^{-1})''(h(\xi))$. 
The necessary details are furnished in Hjort and Fenstad (1993). 
Elementary rules show further that $k_2/k_1^2$ 
also can be written as $-h''(\xi)/h'(\xi)$. 
The special case $h(\xi)=\xi$ has $k_2=0$ and gives back (2.8). 

To analyse the problem associated with (4.3) we can put 
$h(\xi)=\xi^{1/2}$, $d=0$, and use $\xi=1$ in the end. This gives 
$$\lambda_0(c)=\lim_{\eps\arr0}\E\{Q_\eps(c)-Q_\eps(0)\}
	=c^2/\tau^2+((2/3)\gamma/\tau-3/2)c, $$
in terms of standard deviation $\tau$ and skewness $\gamma$.
In the particular $\chi^2_1$ case we find 
$\lambda_0(c)=\half c^2-\sixth c$, with minimum occurring for $c_0=\sixth$. 
We may conclude that the estimator sequence
$\tilda\sigma_N=\{(N-(5/6))^{-1}\sum_{i=1}^N(Y_i-\bar Y_N)^2\}^{1/2}$
can be expected to make fewest errors of the type 
$|\hatt\sigma_N/\sigma-1|\ge\eps$. 
Let us also use (4.4) to exhibit the estimator sequence that 
can be expected to produce the fewest $\eps$-errors on log scale.
So let $Q_\eps^*(c)$ count instances of 
$|\log\hatt\sigma_N^2(c)-\log\sigma^2|\ge\eps$. Then the limit 
can be shown to be $\lambda_0^*(c)=\half c^2-({2\over3}-e^{-1})c$,
with minimum for $c_0={2\over3}-e^{-1}$. So using 
denominator $N-1+c_0$, i.e.
$\sigma_N^*=\{(N-0.695)^{-1}\sum_{i=1}^N(Y_i-\bar Y_N)^2\}^{1/2}$
can be expected to make fewest errors for $\sigma^p$ on the log scale,
for every value of $p$.  

\smallskip
{\csc Remark.} 
Formula (3.2) with $N-1/3$ and the two recently developed 
$N-5/6$ and $N-0.695$ add three new items to the distinguished list 
of denominators in the variance estimation formula 
$\sum_{i=1}^N(Y_i-\bar Y_n)^2/(N-1+c)$ 
reached by various statistical principles. This list includes: 
(i) 
The maximum likelihood estimator, 
normal-based or nonparametric, has $1/N$.
(ii) 
The unbiased estimator for $\sigma^2$ with smallest variance,
under normality or under nonparametric circumstances, 
has $1/(N-1)$. 
(iii) 
If one wants $\E\hatt\sigma_N(c)=\sigma$,
unbiasedness on natural scale, then 
one needs $\E(\chi^2_{n})^{1/2}/(n+c)^{1/2}=1$, 
using $n=N-1$ again, whose approximate solution is $c=-\half$, 
leading to $1/(N-{3\over2})$.
(iv) 
Another possibility is to minimise mean squared error 
$\E_\sigma(\hatt\sigma_N^2-\sigma^2)^2$. 
The best constant is $1/(N+1)$. 
This also gives the best invariant estimator 
under squared error loss for $\sigma^2$. 
(v) 
It appears as natural to solve this problem using natural scale,
i.e.~to minimise $\E_\sigma(\hatt\sigma_N-\sigma)^2$.  
This can again be done uniformly in $\sigma$, and 
the solution is $c=n^2/\{\E(\chi^2_n)^{1/2}\}^2-n$, for which $\half$ is
an approximation. The best denominator is $N-\half$, and gives
at the same time the best invariant estimator under squared error loss
for $\sigma$. 
(vi) 
The median of $\chi^2_n$ is $n-2/3+(4/27)/n$ to a good approximation, 
by the Wilson--Hilferty formula. 
Hence denominator $N-{5\over3}+{4\over 27N}\doteq N-{5\over3}$ 
gives an approximate median-unbiased estimator;
it overshoots as often as it undershoots. 
(vii) 
Since $\E\log(\chi^2_n)=\log2+\psi(n/2)\doteq\log n-({1\over n}+{1\over 3n^2})$
one finds that $\tilda\sigma$ with denominator 
$N-2+{1\over 6N}\doteq N-2$ gives unbiasedness on log scale. 
This is a good property both because the log of the scale is a natural
quantity and because $\tilda\sigma^p$ becomes log-unbiased for $\sigma^p$ 
for every value of $p$. 
(viii) 
The Bayes solution under a vague prior ($\log\sigma$ uniform on the line)
uses $N-3$ if loss is squared error on $\sigma^2$ scale and 
$N-1$ if loss is squared error on $1/\sigma^2$ scale. \square

\bigskip
{\bf 5. Loss functions and Bayes solutions.}

{\sl 5A. Decision theoretic framework.} 
Who can understand his errors (Psalm 19\thinspace:\thinspace12)? 
And who can count them?  
We are essentially working with 
the somewhat non-standard loss function $L_\eps$ 
that for a given sequence of estimates $\{\hatt\theta_n\}$ 
counts 
$$L_\eps\bigl[\theta,\{\hatt\theta_n\colon n\ge1\}\bigr]
	=\sum_{n\ge a(\eps)/\eps^2}
		I\{|\hatt\theta_n-\theta|\ge\eps\}, \eqno(5.1)$$
the number of $\eps$-errors  
among all $n\ge a(\eps)/\eps^2$ cases, for a very small $\eps$. 
The associated risk function is 
$$R_\eps(\theta)=\E_\theta
	L_\eps\bigl[\theta,\{\hatt\theta_n\colon n\ge1\}\bigr]
	=\sum_{n\ge a(\eps)/\eps^2}
       	{\rm Pr}_\theta\{|\hatt\theta_n-\theta|\ge\eps\}. \eqno(5.2)$$
In some cases this series diverges while the difference risk
$R_{1,\eps}(\theta)-R_{2,\eps}(\theta)$ (say) 
might constitute a convergent series. 
The $\lambda_a(c)$ and $\lambda_0(c)$ functions 
found in Section 2 are indeed limits of 
$R_{c,\eps}(\theta)-R_{0,\eps}(\theta)$ as $\eps\rightarrow0$, 
where $R_{c,\eps}$ refers to the risk
function for the $\hatt\theta_n(c)$ sequence. 
One might next obtain average risk with respect to a weight function 
over the parameter space, 
which one may choose or not choose to interpret as 
a prior distribution in the Bayesian fashion,
and finally minimise this expression w.r.t.~the class of estimators
under consideration. 

This was the programme partially carried out in Examples 3A and 4A. 
A more direct and general approach 
which avoids restriction to a given class of estimators 
is possible, as follows, provided the prior is {\it fully} specified. 
Let $\pi(.)$ be the prior distribution and 
$\pi_n(.)=\pi_n(.|x_1,\ldots,x_n)$ 
the posterior distribution at step $n$. Then 
$$\int R_\eps(\theta)\pi(\theta)\,{\rm d}\theta
	=\sum_{n\ge a/\eps^2}\int\cdots\int 
		\pi_n\{|\theta-\hatt\theta_n|\ge\eps\}
      		\,f_n(x_1,\ldots,x_n)\,{\rm d}x_1\cdots {\rm d}x_n,$$
where $f_n$ is the marginal density resulting from having $\pi(.)$ as prior. 
{\it For given $\eps$} the best solution is to minimise each term, 
which means choosing $\hatt\theta_n$ to maximise 
$\pi_n\{|\theta-\hatt\theta_n|\le\eps\}$, for each given $x_1,\ldots,x_n$. 
When $\eps$ goes to zero for fixed $n$ this would mean
using the posterior mode. The balance is more delicate here,
however, where summing over $n$ comes before letting $\eps$ tend to zero. 
Think of $\hatt\theta_n$ as $\theta_n^*+u/\sqrt{n}$, 
where $\theta_n^*$ is the posterior mean for $\theta$. 
Thus one needs to minimise 
$\pi_n\{|\sqrt{n}(\theta-\theta_n^*)-u|\ge\sqrt{s}\}$ w.r.t.~$u$,
where $m=1/\eps^2$ and $s=n/m$ as in earlier sections. 
A full analysis would call for Edgeworth expansions again. 
Observe in particular that the skewness 
of the posterior distribution plays a r\^ole here. 


\smallskip
{\sl 5B. A normal mean.} 
Consider once more the situation of 3A, where it was noted that 
$\theta_n^*$ of (3.1) was best among all linear 
functions of the sample mean, for any prior weight function 
with first moment $\theta_0$ and variance $\tau^2$. 
Suppose now that $\theta$ is given
the normal prior with these parameters.
Then $\pi_n(.)$ is also normal, 
with mean equal to $\theta_n^*$ of (3.1)
and variance $\tau_n^2=\tau^2/(n\tau^2+1)$. The symmetry of the
normal density shows that $\pi_n\{|\theta-\hatt\theta_n|\le\eps\}$
is maximised for $\hatt\theta_n=\theta_n^*$. So (3.1) 
is not only the best solution within the linear class $\hatt\theta_n(c,d)$ 
but the very best of {\it all} estimators, 
under the average $EQ_\eps$ criterion, 
provided the prior used is indeed $N\{\theta_0,\tau^2\}$. 
This even holds for each positive $\eps$ as well as in the limit.

%

\bigskip
{\bf 6. Second order limits in distribution.} 
The starting point for the present work was the first-order 
result (1.1) of Hjort and Fenstad (1992),
relating the distribution of $\eps^2Q_\eps$ to the amount
of time Brownian motion $W(s)$ spends outside the 
$|W(s)|\le s/\sigma$ area. For two estimator sequences 
with the same limiting distribution the typical situation is
that $\eps^2(Q_{1,\eps}-Q_{2,\eps})\arr_p0$ and 
$Q_{1,\eps}/Q_{2,\eps}\arr_p1$, and in the course of the paper 
we have provided limit formulae for {\it the mean value} 
of $Q_{1,\eps}-Q_{2,\eps}$ in different situations. 
This final section briefly presents some 
second order {\it distributional} aspects of this difference,
thus providing additional information about what happens
when $\eps\arr0$.   
In several situations this distribution is in the limit 
related to the amount of time Brownian motion spends along the two 
border lines $\pm s/\sigma$. Such distributions,
for Brownian motion and for general diffusion processes,
are studied in Hjort and Khasminskii (1993). 

To illustrate the typical behaviour, 
let us go back to the situation of Section 2.
Take $c$ and $\xi$ to be positive so that $l(c)>l(0)$ and $r(c)>r(0)$
in (2.3). It follows from (2.4) that 
$$Q_\eps(c)-Q_\eps(0)=\sum_{n\ge a/\eps^2}
	\Bigl[I\bigl\{l(0)\le T_n\le l(c)\bigr\}
	     -I\bigl\{r(0)\le T_n\le r(c)\bigr\}\Bigr]. \eqno(6.1)$$
This is an infinite sum of quite rare 1's and $-1$'s. 
Our previous labours were confined to the expected value
of this difference and therefore stayed within the realm 
of probabilities and Edgeworth expansions. 
Careful analysis of a different type is called for to
find the limit distribution of the (6.1) difference. 
Details are provided partly in Hjort and Khasminskii (1993, Section 6.3) 
and more fully in the technical report Hjort and Fenstad (1993). 
The bottom line is that 
$$\eps\{Q_\eps(c)-Q_\eps(0)\}\arr_d A-B, \eqno(6.2)$$
say, where $A$ and $B$ are two variables concerned 
respectively with how much time Brownian motion spends 
along $s/\sigma$ and $-s/\sigma$. In the case of a fixed positive $a$ 
in (6.1) each of $A$ and $B$ is a mixture of an exponential
and a unit point mass at zero, and the two are dependent. 
In particular $A-B$ can then be zero with positive probability. 
In the case of a shrinking $a=a(\eps)$, say $a=\eps$ in (6.1), 
both $A$ and $B$ are exponentially distributed with mean $c\xi$ 
and with intercorrelation $-{1\over3}$.

		
\bigskip
\parindent0pt\parskip3pt\baselineskip12pt

\centerline{\bf References}

\medskip
\ref{%
Barndorff-Nielsen, O.E.~and Cox, D.R. (1989).
{\sl Asymptotic Techniques for Use in Statistics.}
Chapman and Hall, London.} 

\ref{%
Berkson, J. (1980).
Minimum chi-square, not maximum likelihood!
(with discussion). 
{\sl Ann.~Statist.}\allowbreak~{\bf 8}, 457--487.} 



\ref{%
Ghosh, J.K.~and Subramanyam, K. (1974).
Second order efficiency of maximum likelihood estimators. 
{\sl Sankhy\=a, Ser.~A} {\bf 36}, 325--258.} 

\ref{%
Hodges, J.~and Lehmann, E.L. (1970).
Deficiency.
{\sl Ann.~Math.~Statist.}~{\bf 41}, 783--801.}

\ref{%
Hjort, N.L.~and Fenstad, G. (1992).
On the last time and the number of times an estimator
is more than $\eps$ from its target value.
{\sl Ann.~Statist.}~{\bf 20}, 469--489.} 

\ref{%
Hjort, N.L.~and Fenstad, G. (1993).
Second order asymptotics for the number of times an estimator 
is more than $\eps$ from its target value.
Statistical Research Report, Department of Mathematics, 
University of Oslo.} 

\ref{%
Hjort, N.L.~and Khasminskii, R.Z. (1993).
On the time a diffusion process spends along a line.
{\sl Stoch.~Proc.~and their Appl.}~{\bf 47}, 229--247.} 

\ref{%
Kolassa, J.E.~and McCullagh, P. (1990).
Edgeworth series for lattice distributions.
{\sl Ann.\allowbreak~Statist.}~{\bf 18}, 981--985.} 

\ref{%
Lehmann, E.L. (1983).
{\sl Theory of Point Estimation.}
Wiley, New York.} 

\ref{%
Pfanzagl, J. (1973). 
Asymptotic expansions related to minimum contrast estimators.
{\sl Ann. Statist.}\allowbreak~{\bf 1}, 993--1026. }

\ref{%
Rao, C.R. (1962).
Efficient estimates and optimum inference procedures in large samples
(with discussion). 
{\sl J.~Roy.~Statist.~Soc., Ser.~B} {\bf 24}, 46--72. }

\bye